\documentclass[12pt]{article}
\usepackage{amsmath}
\usepackage{amssymb}
\usepackage{amscd}
\newtheorem{theorem}{Theorem}
\newtheorem{definition}{Definition}
\newtheorem{proposition}{Proposition}
\newtheorem{corollary}{Corollary}
\newtheorem{lemma}{Lemma}
\newtheorem{example}{Example}
\newtheorem{remark}{Remark}
\newenvironment{proof}{\begin{trivlist} \item[]{\it Proof---}}
{\par\hfill $\square$\end{trivlist}}
\renewcommand{\P}{\mathbb{P}}

\newcommand{\C}{\mathbb{C}}

\renewcommand{\H}{{\cal H}}

\newcommand{\supp}{{\rm supp}}

\newcommand{\E}{{\cal E}}
\newcommand{\F}{{\cal F}}

\renewcommand{\d}{{\rm d}}
\newcommand{\ddc}{{\rm dd}^{\rm c}}
\newcommand{\ddbar}{\partial\overline\partial}
\newcommand{\dist}{{\rm dist}}
\newcommand{\Loneloc}{{\rm L}^1_{\rm loc}}

\newcommand{\Ball}{{\rm B}}
\newcommand{\hull}{{\rm hull}}

\newcommand{\see}{{\it see }}
\newcommand{\ie}{{\it i.e. }}
\newcommand{\link}{{\rm link}}
\title{Polynomial hulls and positive currents}
\author{Tien-Cuong Dinh}
\begin{document}
\maketitle
\begin{abstract} Let $T$ be a positive plurisubharmonic current of bidimension 
$(p,p)$ and let $\delta>0$. Assume that the Lelong number of $T$ satisfies 
$\nu(T,a)\geq \delta$ on a dense subset of $\supp(T)$ (rectifiable
currents satisfy this condition). 
Then $T=\varphi[X]$, where 
$X$ is a complex subvariety of pure dimension $p$ and $\varphi$ is a weakly 
plurisubharmonic function on $X$. We have an 
analogous result for plurisuperharmonic currents.
We also introduce and study a notion of polynomial $p$-hull.
\end{abstract}
\section{Introduction}
Let $V$ be a complex manifold of dimension $n$. 
A current $T$ of bidimension $(p,p)$
or bidegree $(n-p,n-p)$ in $V$ is called 
{\it plurisubharmonic} if $\ddc T$ is a
positive current, is called 
{\it plurisuperharmonic} if $-\ddc T$ is a
positive current
and is called {\it pluriharmonic} if $\ddc T=0$. 
In particular, every closed current is pluriharmonic.
In this paper, we study some properties of positive
plurisubharmonic and plurisuperharmonic currents. 
\par
Positive plurisubharmonic and plurisuperharmonic currents are very 
useful in many mathematical problems. In 
\cite{Garnett}, Garnett proved that every laminated 
closed set supports a positive
pluriharmonic currents (see also \cite{BerndtssonSibony}). Duval and Sibony 
used positive plurisuperharmonic currents to describe the polynomial hull of compact
sets in $\C^n$ \cite{DuvalSibony}. Gauduchon,
Harvey, Lawson, Michelson, Alessandrini, Bassanelli, etc. 
studied non-K\"ahler geometry using (smooth or not) pluriharmonic 
currents \cite{AB1}. 
\par
Many important properties of these classes of currents have been considered in the
papers of 
Skoda \cite{Skoda2}, Sibony, Berndtsson, Forn\ae ss \cite{Sibony},
\cite{BerndtssonSibony},
\cite{FornaessSibony}, Alessandrini and Bassanelli 
\cite{AB2,AB3,Bassanelli}, etc. 
In particular, we know that if $T$ is a positive plurisubharmonic
current, one can define a density $\nu(T,a)$ of $T$ at every point $a$. 
This density is called Lelong number of $T$ at $a$. 
\par
We will prove that if 
the level set $\{\nu(T,a)\geq\delta\}$ is dense 
in the support of $T$ for a suitable 
$\delta>0$, then the support
$X$ of $T$ is a complex subvariety of pure dimension $p$ 
and $T=\varphi[X]$, where $\varphi$ 
is a weakly plurisubharmonic function on $X$. 
In particular, every rectifiable 
positive plurisubharmonic current is closed. Moreover, it has the form 
$c[X]$, where $X$ is a complex
subvariety and $c$ is essentially equal to a positive 
integer in each component of $X$.  
We will also prove that
if $T$ is positive pluriharmonic (resp. plurisuperharmonic)
and its support has locally finite 
$2p$-dimensional Hausdorff measure then the support
$X$ of $T$ is a complex subvariety of pure dimension $p$ 
and $T=\varphi[X]$, where $\varphi$ 
is a weakly pluriharmonic (resp. plurisuperharmonic) function on $X$. 
\par
For closed positive currents, an analogous result 
has been proved by King \cite{King}
(see also \cite{Skoda, Siu}). The proof of King does not work in the 
case of plurisubharmonic and plurisuperharmonic currents. 
Our proof relies on some recent progress on the polynomial hull of 
finite length compact sets \cite{Dinh1,Dinh2,Lawrence2}. We also use 
a theorem of Shiffman on separately
holomorphic functions and some results of 
Alessandrini-Bassanelli on the Lelong number.
\par
In order to study the currents of bidimension $(p,p)$ we will introduce 
the notion of polynomial $p$-hull which is the usual polynomial hull when $p=1$. 
We will extend the Wermer theorem for $p$-hull.
\section{Polynomial hulls and pseudoconcave sets}
Let $\Gamma\subset \C^n$ be a compact set.
{\it The polynomial hull} of $\Gamma$ is the compact set $\widehat \Gamma$
defined by the following formula
$$\widehat \Gamma:=\big\{x\in \C^n,\ |P(x)|\leq \max_{z\in\Gamma} |P(z)| 
\mbox{ for any polynomial } P \big\}.$$ 
We now introduce the notions of $p$-pseudoconcavity and polynomial $p$-hull.
\par  
Let $V$ be a complex manifold of dimension $n\geq 2$ and $X$ a closed subset of $V$. 
Let $1\leq p\leq n-1$ be an integer. 
\begin{definition} \rm
We say that $X$ is {\it $p$-pseudoconcave
subset of $V$}
if for every open set $U\subset\subset V$ and every holomorphic map $f$ from a 
neighbourhood of $\overline U$ into
$\C^p$ we have $f(X\cap U)\subset \C^p\setminus \Omega$ 
where $\Omega$ 
is the unbounded component of $\C^p\setminus f(X\cap bU)$. 
\end{definition}
This means $X$ has no ``locally peak point'' for holomorphic 
maps into $\C^p$. In particular, $X$ satisfies the maximum principle, i.e. 
holomorphic functions admit no local strict maximum modul on $X$. 
Observe that if $V$ is a submanifold of another complex manifold $V'$, then $X$ is 
$p$-pseudoconcave in $V$ if and only if $X$ is $p$-pseudoconcave in $V'$.
\par
By the argument principle, every complex subvariety 
of pure dimension $p$ of $V$ is 
$p$-pseudoconcave. It is clear that 
if $X$ is $p$-pseudoconcave, the $2p$-dimensional Hausdorff measure of $X$ 
is strictly positive. 
If $1\leq p\leq q\leq n-1$ and $X$ is $q$-pseudoconcave then $X$ is 
$p$-pseudoconcave. If $g:V\longrightarrow \C^{p-k}$ is a holomorphic map
 and $X\subset V$ is $p$-pseudoconcave then $g^{-1}(x)\cap X$ is 
$k$-pseudoconcave for every $x\in \C^{p-k}$. We have the following proposition.
\begin{proposition} Let $T$ be a positive plurisuperharmonic current of bidimension 
$(p,p)$ on a complex manifold $V$. Then the support $\supp(T)$ of 
$T$ is $p$-pseudoconcave in $V$.
\end{proposition}
\begin{proof} Since the problem is local we can assume that $V$ is a ball of $\C^n$.
Assume that $X:=\supp(T)$ is not $p$-pseudoconcave. Then there exist
an open set $U\subset\subset V$ 
and a holomorphic map $f:V'\longrightarrow \C^p$ such that
$f(X\cap U)\not \subset \C^p\setminus \Omega$, where $V'$ 
is a neighbourhood of $\overline U$ and $\Omega$ is 
the unbounded component of $\C^p\setminus f(X\cap bU)$.
\par
Let $\Phi:V'\longrightarrow \C^p\times V'$ the holomorphic map
given by $\Phi(z):=(f(z),z)$. Choose a bounded domain $U'$ in $\C^p\times V'$ such that 
$U'\cap\Phi(V')=\Phi(U)$.
Set $X':=\Phi(X\cap V')$ and $T':=1_{U'}\Phi_*(T)$. Then $T'$ is a 
positive plurisuperharmonic
current in $U'$. 
Let $\pi:\C^{n+p}\longrightarrow \C^p$ be 
the linear projection on the first $p$ coordinates. We have 
$\pi(\overline X'\cap bU')=f(X\cap bU)$ and 
$\pi(X'\cap U')= f(X\cap U)\not \subset \C^p\setminus\Omega$. 
The open set $\Omega$ is
also the unbounded component of $\C^p\setminus 
\pi(\overline X'\cap bU')$.
\par
Therefore $\pi_*(T')$ defines a positive plurisuperharmonic
current in $\Omega$ which is null in $\C^p\setminus\pi(\overline X')$.
Hence there is a positive plurisuperharmonic
function $\psi$ on $\Omega$, 
null on $\C^p\setminus\pi(\overline X')$
such that $\pi_*(T')=\psi[\Omega]$ in $\Omega$. By plurisuperharmonicity, this function is 
identically null.
Fix a small open $W\subset \C^{n+p}$ such that 
$W\cap X'\not=\emptyset$ and $\pi(\overline W)\subset \Omega$. 
Set $\Psi:=(i\d z_1\wedge \d \overline z_1+\cdots+i\d z_p\wedge \d \overline z_p)^p$.
Then the positive measure
$T'\wedge \pi^*(\Psi)$ is null in $W$ since 
its push-forward by $\pi$ is null. This still holds for every small 
linear pertubation $\pi_\epsilon$ of $\pi$. On the other hand, we can construct
a strictly positive $(p,p)$-form $\psi$ of $\C^{n+p}$ as a linear combinaison of 
$\pi_\epsilon^*(\Psi)$. Then we have $\langle T', 1_W \psi\rangle =0$. 
Since $T'$ is positive, $T'=0$ on $W$. This is a contradiction.
\end{proof}
\begin{remark}\rm
We can also prove 
Proposition 1 by using the Oka inequality of Forn\ae ss and Sibony 
\cite{FornaessSibony}
\end{remark}
\begin{proposition} Let $\Gamma$ be a compact subset of $\C^n$. 
Denote by $\F$ the family of $p$-pseudoconcave subsets of $\C^n\setminus\Gamma$ 
which are bounded in $\C^n$. Then the union $\Sigma$ of elements of $\F$ is also belong to
$\F$ (then $\Sigma$ is the biggest element of $\F$). 
\end{proposition}
\begin{proof} 
Let $S_1$, $S_2$, $\ldots$ be elements of $\F$ 
such that $S_i\subset S_{i+1}$. We show that $S:=\overline{\cup S_i}$ belongs to $\F$.
By maximum principle, $S_i\subset\widehat\Gamma$. Then $S$ is bounded in $\C^n$.
\par
Assume that $S$ is does not belong to $\F$. Then there are a point $p\in S$, a 
neighbourhood 
$U\subset\subset\C^n\setminus \Gamma$ of $p$ and a holomorphic map $f$ from a 
neighbourhood 
of $\overline U$ into $\C^p$ such that $f(p)$ belongs to the unbounded component $\Omega$ 
of $\C^p\setminus f(S\cap bU)$. Fix small neighbourhoods $V$ of $p$ and $W$ of
$S\cap bU$ such that $f(\overline V)$ is included in the unbounded component of 
$\C^p\setminus f(\overline W)$. 
\par
For $i$ big enough, we have $S_i\cap bU\subset W$ and $S_i\cap V\not=\emptyset$. 
Moreover, if $p$ is a point in 
$S_i\cap V$, $f(p)$ belongs to the unbounded component of 
$\C^p\setminus f(S_i\cap bU)$. This contradicts the $p$-pseudoconcavity. 
Thus $S\in\F$.
\par
Now, by Zorn axiom, we can choose an element $\Sigma'$ of 
$\F$ which is maximal for the
inclusion. If $\Sigma'\not=\Sigma$ then there is an element $S\in\F$ such that 
$S\not\subset \Sigma'$. The set $S\cup\Sigma'$ is $p$-pseudoconcave in 
$\C^n\setminus\Gamma$ and it belongs to $\F$. This contradicts 
the maximality of $\Sigma'$.
\end{proof}
\begin{definition} \rm
Let $\Gamma$ and $\Sigma$ be as in Proposition 2. We say that 
$\Sigma\cup\Gamma$ is the 
{\it polynomial $p$-hull} of $\Gamma$ and we denote it by $\hull(\Gamma,p)$.
\end{definition}
The following proposition shows that we obtain the usual polynomial hull when $p=1$.
\begin{proposition} We have $\hull(\Gamma,1)=\widehat\Gamma$ for every compact 
subset $\Gamma$ of $\C^n$.
\end{proposition}
\begin{proof} By maximum principle, we have 
$\hull(\Gamma,1)\subset \widehat\Gamma$. Now let 
$z\in\widehat\Gamma\setminus\Gamma$, we will prove that $z\in\hull(\Gamma,1)$. 
By Duval-Sibony theorem \cite{DuvalSibony}, there are
a positive current $T$ of bidimension $(1,1)$ 
with compact support in $\C^n$ and a measure $\mu$ with support
in $\Gamma$ such that $\ddc T=\mu-\delta_z$ where $\delta_z$ 
is the Dirac mass at $z$.
By Proposition 1, $\supp(T)$ is $1$-pseudoconcave in $\C^n\setminus\Gamma$. Thus 
$\supp(T)\subset \hull(\Gamma,1)$ and $z\in \hull(\Gamma,1)$.
\end{proof}
\section{Hulls of finite Hausdorff measure sets}
In this section, we first 
recall some results on the polynomial hull of finite length
compact sets and then we will extend them for polynomial $p$-hull.
\par
Denote by $\H^k$ the Hausdorff measure of dimension $k$. A compact set 
$\Gamma\subset \C^n$ is called {\it geometrically $k$-rectifiable} 
if $\H^k(\Gamma)$ is finite and the geometric tangent cone of $\Gamma$ 
is a real space of dimension $k$ at $\H^k$-almost every point in $\Gamma$.
If a compact set $\Gamma$ is geometrically $k$-rectifiable, it is
{\it $(\H^k,k)$-rectifiable} \cite[p.208]{Mattila}, 
i.e. there exist ${\cal C}^1$ 
manifolds $V_1$, $V_2$, $\ldots$ such that 
$\H^k(\Gamma\setminus \cup V_m)=0$
\cite[3.1.16, 3.2.18]{Federer}.
Now, suppose that $V_1$, $V_2$, $\ldots$ are ${\cal C}^1$ 
oriented manifolds of dimension $k$
in $\C^n$, $K_i\subset V_i$ are compact sets and $n_1$, $n_2$, 
$\ldots$ are integers such that
$\sum |n_i|\H^k(K_i)<+\infty$. Then we can define 
a current $S$ of dimension $k$ by
$$\langle S,\psi \rangle := \sum n_i \int_{K_i} \psi$$
for any test form $\psi$ of degree $k$ having compact support. 
Such a current is called
{\it rectifiable current}. 
\par  
We have the following theorem.
\begin{theorem}
Let $\Gamma$ be a compact subset of $\C^n$. Assume that $\Gamma$ is 
geometrically $(2p-1)$-rectifiable with $1\leq p\leq n-1$. 
Then $\hull(\Gamma,p)\setminus \Gamma$ is a
complex subvariety of pure dimension $p$ (possibly empty) 
of $\C^n\setminus \Gamma$. Moreover, 
$\hull(\Gamma,p)\setminus \Gamma$ has finite $2p$-dimensional Hausdorff measure 
and the boundary of the integration current 
$[\hull(\Gamma,p)\setminus \Gamma]$ 
is a rectifiable current of dimension $2p-1$ and has multiplicity $0$ or $1$ 
$\H^{2p-1}$-almost everywhere on $\Gamma$.
\end{theorem}
The last property of the current 
$[\hull(\Gamma,p)\setminus \Gamma]$ is called 
{\it Stokes formula}. In the case $p=1$, this theorem is proved in 
\cite{Dinh1,Dinh2, Lawrence2} and it generalizes the
results of Wermer \cite{Wermer}, Bishop \cite{Bishop}, Stolzenberg 
\cite{Stolzenberg}, Alexander \cite{Alexander1} and 
Harvey-Lawson \cite{HarveyLawson1} ({\it see} also \cite{DolbeaultHenkin}). 
\par
In order to prove Theorem 1, we will slice $\Gamma$ by complex planes 
of dimension $n-p+1$
and we will apply the known result in the case $p=1$. We will also use a theorem of
Shiffman on separatly holomorphic functions.
\par
Set $X:=\hull(\Gamma,p)\setminus\Gamma$. Observe that if $L$ is a linear 
complex $(n-p+1)$-plane
then $X\cap L$ is $1$-pseudoconcave in $L\setminus\Gamma$. 
Thus $X\cap L$ is included in
the polynomial hull of $\Gamma\cap L$.
Moreover, by Sard theorem, for almost every
$L$ the intersection $\Gamma\cap L$ is geometrically $1$-rectifiable (\see also 
\cite{Dinh1}). Therefore
$\hull(\Gamma\cap L,1)\setminus\Gamma$ is a complex subvariety of 
pure dimension $1$ of 
$L\setminus\Gamma$. We need the following lemma.
\begin{lemma} Let $V$ be a complex manifold of dimension $n\geq 2$ and $X\subset V$ a 
$p$-pseudoconcave subset. If $X$ is included in 
a complex subvariety $\Sigma$ of pure 
dimension $p$ of $V$, then $X$ is itself a complex 
subvariety of pure dimension $p$ of $V$.  
\end{lemma}
\begin{proof} Recall that every non empty $p$-pseudoconcave 
set has positive $\H^{2p}$ 
measure. Assume that $X$ is not a complex 
subvariety of pure dimension $p$ of $V$. 
Then there is a point $a\in X$ such that 
$a$ is a regular point of $\Sigma$ and $a\in \overline{\Sigma\setminus X}$. 
Let $V'\subset V$ 
be a small ball of centre $a$ satisfying 
$V'\cap \Sigma\setminus X\not=\emptyset$. Since $V'$ 
is small, we can choose a projection $\pi:V'\longrightarrow \C^p$ such that 
$\pi$ is injective on $\Sigma\cap V'$. Then 
$\pi(a)$ is belong to the unbounded component
of $\C^p\setminus\pi(X\cap bW)$ for every neighbourhood 
$W\subset\subset V'$ of $a$. This is impossible.
\end{proof}
By the latter lemma, for almost every $L$ (such that $\Gamma\cap L$ is geometrically
$1$-rectifiable), $X\cap L$ is a complex subvariety of pure dimension 1 of 
$L\setminus \Gamma$. Let $E$ be a complex 
$(n-p)$-plane which does not meet $\Gamma$. 
By Sard theorem, for almost every $L$ passing through $E$, 
the intersecrtion $\Gamma\cap L$
is geometrically $1$-rectifiable. We deduce that 
$X\cap E=(X\cap L)\cap E$ is a finite set.
We now use a Shiffman theorem on separately 
holomorphic functions in order to 
complete the proof of Theorem 2. 
\par
Let $a$ be a point of $X$. We show that $X$ is a complex 
variety of pure dimension $p$ 
in a neighbourhood of $a$. 
Indeed, choose a coordinates system such that 
$\Pi(a)\not\in \Pi(\Gamma)$,
where $\Pi:\C^n\longrightarrow \C^p$ is the linear projection given by
$\Pi(z):=(z_1,\ldots, z_p)$. Set $\pi_i:\C^n\longrightarrow \C^{p-1}$,
$\pi_i(z):=(z_1,\ldots,z_{i-1},z_{i+1},\ldots,z_p)$. 
Using a linear change of coordinates, 
we may suppose without loss of generality that for $\H^{2p-2}$-almost 
every $x\in\C^{p-1}$ and for every $1\leq i\leq n$, 
$X\cap\pi_i^{-1}(x)$ is a complex subvariety of pure dimension $1$ of 
$\pi_i^{-1}(x)\setminus \Gamma$. 
\par
Fix a small 
open neighbourhood $V$ of $\Pi(a)$
such that $\Pi^{-1}(V)\cap\Gamma=\emptyset$.
Put $z=(z',z'')$, 
$z':=(z_1,\ldots,z_p)$ and $z'':=(z_{p+1},\ldots,z_n)$. 
Let $g$ be a function which is defined in a subset $V'$ 
of total measure of $V$. 
This function is called {\it separately holomorphic} if 
for every $1\leq i\leq p$ 
and $\H^{2p-2}$-almost every $x\in \C^{p-1}$ the restriction of $g$ on 
$V'\cap \pi_i^{-1}(x)$ can be extended to a holomorphic 
function on $V\cap \pi_i^{-1}(x)$.
For any $k\geq 0$ and 
$p+1\leq m\leq n$, we define the following measurable function on $V$ 
$$g_{m,k}(z'):=\sum_{z\in\Pi^{-1}(z')\cap X} z_m^k.$$
The function
$$g_0(z'):= g_{m,0}(z)=\#(\Pi^{-1}(z')\cap X)$$
is independent on $m$ and takes only positive integer values. 
\par
By the choice of coordinates, $g_{m,k}(z')$ is $\H^{2p}$-almost everywhere equal to a 
separately holomorphic function. Thanks to Shiffman Theorem 
\cite{Shiffman}, $g_{m,k}(z')$ is equal
$\H^{2p}$-almost everywhere to a function $\tilde g_{m,k}(z')$ 
which is holomorphic on $V$.
In particular, $\tilde g_{m,0}(z')$ is equal to an integer $r$ 
which does not depend on $m$ and $z'$.    
\par
Now, consider the following equation system:
$$\sum_{i=1}^r \big[z_m^{(i)}]^k = \tilde g_{m,k}(z') 
\mbox{ \ \ with } p+1\leq m\leq n 
\mbox{ and } 1\leq k\leq r.$$
The set of points $\big( z', z_{p+1}^{(i_{p+1})}, \ldots, z_n^{(i_n)} \big)$ 
given by solutions of the system above is a 
complex subvariety $\Sigma$ of pure dimension $p$ of $\Pi^{-1}(V)$. 
We have $X\cap\Pi^{-1}(x)\subset\Sigma$ for $\H^{2p}$-almost every $x\in V$. 
We show that
$X\cap\Pi^{-1}(V)\subset \Sigma$. 
Asssume this is not true. Choose
a point $b\in X\cap \Pi^{-1}(V)\setminus \Sigma$ and a neighbourhood 
$W\subset\subset\Pi^{-1}(V)\setminus \Sigma$ of $b$. 
The set
$\Pi(X\cap \overline W)$ has zero $\H^{2p}$ measure.
But for almost every complex $(n-p+1)$-plane 
$L$ passing through $\Pi^{-1}(\Pi(b))$, 
the intersection $X\cap L$ is a complex subvariety of pure dimension $1$ of $L\setminus
\Gamma$ which contains $b$. This implies that 
the set $\Pi(X\cap \overline W)$ 
has positive $\H^{2p}$ measure.
We reach a contradiction.
\par
Now we have already know that $X\cap\Pi^{-1}(V)\subset\Sigma$. 
By Lemma 1, $X\cap \Pi^{-1}(V)$ is a complex 
subvariety of pure dimension $p$
of $\Pi^{-1}(V)$. 
Then $X$ is a complex subvariety of pure dimension $p$ of $\C^n\setminus 
\Gamma$. The rest of the proof follows along the sames lines as the proof of 
Stokes formula given in \cite{Lawrence1, Dinh1, Dinh2}.
\par
\hfill $\square$
\\
When $\Gamma$ is included in the boundary of a smooth convex domain, its rectifiability 
in Theorem 1 is non necessary. We have the following result.
\begin{theorem} Let $D$ be a strictly convex 
domain with ${\cal C}^2$ boundary of $\C^n$. Let $\Gamma$ be a compact 
set in the boundary of $D$. 
Assume that the $(2p-1)$-dimensional Hausdorff measure  
of $\Gamma$ is finite. Then $\hull(\Gamma,p)\setminus \Gamma$ is a 
complex subvariety of pure dimension $p$ (possibly empty) 
of $\C^n\setminus \Gamma$. Moreover, 
$\hull(\Gamma,p)\setminus \Gamma$ has finite $2p$-dimensional Hausdorff measure 
and the boundary of the integration current 
$[\hull(\Gamma,p)\setminus \Gamma]$ 
is a rectifiable current of dimension $2p-1$ and has multiplicity $0$ or $1$ 
$\H^{2p-1}$-almost everywhere on $\Gamma$.
\end{theorem}
In the case $p=1$, this result is proved by Lawrence \cite{Lawrence2} as a corollary 
of Theorem 1. In order to apply Theorem 1, Lawrence show that if a compact set
$\Gamma'\subset\Gamma$ is minimal
with respect to the property that its polynomial
hull contains some point $w\in D$, then it is 
geometrically $1$-rectifiable. Using this result of Lawrence, the proof of Theorem 1
remains valid for Theorem 2.
\begin{remark} \rm
Under the hypothesis of Theorems 1 and 2, 
by Lemma 1, every $p$-pseudoconcave 
subset of $\C^n\setminus \Gamma$ which is bounded in $\C^n$, 
is a complex subvariety of 
pure dimension $p$ of $\C^n\setminus\Gamma$ and is included in $X$. 
By Proposition 1, every positive 
current $T$ with compact support satisfying $\ddc T\leq 0$ in 
$\C^n\setminus\Gamma$ has
the form $\varphi[X]$, where $\varphi$ is a weakly plurisuperharmonic 
function in $X$ 
(\see the next section). 
\end{remark}
\begin{corollary} Let $V$ be a complex manifold of dimension $n\geq 2$ and $X$ a 
$p$-pseudoconcave subset of
$V$. Let $K$ be a compact subset of $V$ which admits a Stein neighbourhood.
Assume that the $2p$-dimensional Hausdorff measure  
of $X\setminus K$ is locally 
finite in $V\setminus K$. Then $X$ is a complex subvariety of 
pure dimension $p$ of $V$. 
\end{corollary}
\begin{proof} Fix a point $a\in X\setminus K$ and $\Ball\subset \subset V\setminus K$ an 
open ball contening $a$ such that $X\cap b\Ball$ has finite $\H^{2p-1}$ measure. 
Then $X\cap \Ball$ is $p$-pseudoconcave in $\Ball$. We can consider $\Ball$ as 
an open ball of $\C^n$. Then $X\cap \Ball\subset \hull(X\cap b\Ball,p)$. 
By Remark 2,
$X\cap \Ball$ is a complex subvariety of pure dimension $p$ of $\Ball$. 
This implies that $X\setminus K$ is
a complex subvariety of pure dimension $p$ of $V\setminus K$.
\par
Since we can replace $V$ by a Stein neighbourhood of $K$, suppose without loss of
generality that $V$ is a submanifold of $\C^N$. 
Observe that $X$ is also $p$-pseudoconcave
in $\C^N$. Let $\Ball'$ be a ball containing the compact $K$ 
such that $X\cap \Ball'$ has 
finite $\H^{2p-1}$ measure. Then we obtain 
$X\cap \Ball'\subset \hull(X\cap b\Ball',p)$. 
By Remark 2, $X\cap \Ball'$ is a complex subvariety of pure 
dimension $p$ of $\Ball'$. This
completes the proof.
\end{proof}
\section{Positive currents}
Let $V$ be a complex manifold of dimension $n$. 
Let $X$ be a complex subvariety of pure 
dimension $p$ of $V$. 
A function $\varphi\in \Loneloc(X)$ is called
{\it weakly plurisubharmonic} if it is essentially bounded from above and 
$\ddc(\varphi[X])$ is a positive current.
In the regular part of $X$, $\varphi$ is equal almost everywhere to an 
usual plurisubharmonic function. We say that 
$\varphi$ is 
{\it weakly plurisuperharmonic} if $-\varphi$
is weakly plurisubharmonic 
and that $\varphi$ is 
{\it weakly pluriharmonic} if both $\varphi$ and $-\varphi$
are weakly plurisubharmonic. 
\par
If $\varphi$ is a positive weakly 
plurisubharmonic function on $X$, we can define 
a positive plurisubharmonic current $T:=\varphi[X]$ as follows:
$$\langle T,\psi\rangle := \int_X \varphi\psi$$
for every test form $\psi$ of bidegree $(p,p)$ compactly supported in $V$.
According to Bassanelli \cite[1.24, 4.10]{Bassanelli}, 
every positive plurisubharmonic current with 
support in $X$ can be obtained by this way. This is also true for positive 
plurisuperharmonic and pluriharmonic currents.
\par
Now, let $T$ be a positive plurisubharmonic current of bidimension $(p,p)$ 
of $V$. 
Thanks to a Jensen type formula \cite{AB3,Skoda2,Demailly}, 
if $V$ is an open set of $\C^n$
which contains a ball $\Ball(a,r_0)$, the function
$$\frac{1}{2^p\pi^p r^{2p}}\langle T,1_{\Ball(a,r)}\omega^p\rangle, 
\mbox{ for } 0\leq r <r_0,$$
is non-decreasing, where $\omega:=i \d z_1\wedge \d\overline z_1+\cdots + 
i \d z_n\wedge \d\overline z_n$. Therefore, we can define
$$\nu(T,a):=\lim_{r\rightarrow 0} 
\frac{1}{2^p\pi^p r^{2p}}\langle T,1_{\Ball(a,r)}\omega^p\rangle.$$
This limit is called {\it Lelong number} of $T$ at the point $a$. 
In \cite{AB3},
Alessandrini and Bassanelli proved that 
$\nu(T,a)$ does not depend on coordinates. Thus, the Lelong number 
is well defined for every manifold $V$.
\par
We have the following theorem.  
\begin{theorem} Let $T$ be a positive plurisubharmonic current 
of bidimension $(p,p)$ in a complex manifold $V$ of dimension $n$. 
Assume that there exists a real number $\delta>0$ such that the level set
$\{z\in V,\ \nu(T,z)\geq \delta\}$ is dense in the support 
$\supp(T)$ of $T$. Then $\supp(T)$ is a complex subvariety of pure 
dimension $p$ of $V$ and there exists a weakly plurisubharmonic 
function $\varphi$ on $X:=\supp(T)$ such that $T=\varphi[X]$.  
\end{theorem}
\begin{remark} \rm
In the case of closed positive currents, an analogous result is proved 
by King \cite{King}. 
\end{remark}
It is sufficient to prove that the support $X=\supp(T)$ of $T$ is a complex
subvariety of pure dimension $p$ of $V$. Since the problem is local, 
we can suppose 
that $V$ is a ball of $\C^n$. By Corollary 1, we have to prove 
that $X$ has locally 
finite $\H^{2p}$ measure and $X$ is $p$-pseudoconcave. Consider now  
the {\it trace measure} $\sigma:= T\wedge\omega^p$ of $T$.
\begin{proposition} 
Let $T$ be a plurisubharmonic current of bidimension $(p,p)$ 
in a complex manifold
$V$. Then for every $\delta>0$ the level set 
$\{\nu(T,a)\geq \delta\}$ is closed
and has locally finite $\H^{2p}$ measure. Moreover, 
under the hypothesis of Theorem 3,
we have $\nu(T,a)\geq\delta$ for every $a\in X$ and $X$ has locally 
finite $\H^{2p}$ measure.
\end{proposition}
\begin{proof} We may suppose without loss of generality 
that $V$ is an open subset of $\C^n$. 
Set $Y:= \{\nu(T,a)\geq \delta\}$. Let $(a_n)\subset Y$ be a sequence which 
converges to a point $a\in V$.
Fix an $r>0$ such that $\Ball(a,r)\subset V$. 
We have 
$$\sigma(\Ball(a,r))\geq\sigma(\Ball(a_n,r-|a_n-a|))\geq 2^p\pi^p\delta 
(r-|a_n-a|)^{2p},$$
which implies that
$$\sigma(\Ball(a,r))\geq 2^p\pi^p\delta r^{2p}$$
and
$$\nu(T,a)=\lim_{r\rightarrow 0}\frac{\sigma(\Ball(a,r))}{2^p\pi^p
  r^{2p}}\geq \delta.$$
Therefore, $a\in Y$ and hence $Y$ is closed.
The last inequality also implies that $Y$ has locally finite $\H^{2p}$ measure 
\cite[2.10.19(3)]{Federer}.
\end{proof}
\begin{lemma} Let $T$ be a positive current of bidimension $(p,p)$ 
in $\C^n$ with compact support. Let $\pi:\C^n\longrightarrow 
\C^p$ be the linear 
projection on the first $p$ coordinates, $a\in \C^n$ and $b:=\pi(a)$.
Assume that $T$ is plurisubharmonic on $\pi^{-1}(\Ball(b,r))$ 
for some $r>0$.
Then $\langle T,1_{\pi^{-1}(\Ball(b,r))}\Psi\rangle 
\geq 2^p\pi^p r^{2p} \nu(T,a)$,
where $\Psi:=(i\d z_1\wedge \d \overline z_1+ \cdots + 
i\d z_p\wedge \d\overline z_p)^p$. In particular, 
the Lelong number of $\pi_*(T)$
at $b$ is greater or equal to $\nu(T,a)$.
\end{lemma}
\begin{proof} We may suppose without loss of generality $a=0$ and $b=0$.
Set $z':=(z_1,\ldots,z_p)$, $z'':=(z_{p+1},\ldots, z_n)$, 
$\omega:=i\d z_1 \wedge \d 
\overline z_1+ \cdots + i\d z_n\wedge \d\overline z_n$ and
fix an $\epsilon>0$. Let $\Phi:\C^n\longrightarrow \C^n$ 
be the linear map given 
by $\Phi_\epsilon(z):=(z',\epsilon z'')$ and set $T':=(\Phi_\epsilon)_*T$. 
Since the Lelong number
is independant on coordinates \cite{AB3}, we have
$\nu(T',0)=\nu(T,0)$. 
By Jensen
type formula we have 
$$\langle T',1_{\pi^{-1}(\Ball(0,r))}\omega^p\rangle 
\geq 2^p\pi^p r^{2p} \nu(T',0)$$ 
since $\pi^{-1}(\Ball(0,r))$ contains the ball $\Ball(0,r)$ of $\C^n$.
From this inequality, it follows that
$$\langle T,1_{\pi^{-1}(\Ball(0,r))}\Phi_\epsilon^*(\omega^p)\rangle 
\geq 2^p\pi^p r^{2p} \nu(T,0).$$ 
On the other hand, when 
$\epsilon\rightarrow 0$, $\Phi_\epsilon^*(\omega^p)$ converges uniformly to 
$\Psi$. Thus,
$$\langle T,1_{\pi^{-1}(\Ball(0,r))}\Psi\rangle 
\geq 2^p\pi^p r^{2p} \nu(T,0).$$ 
This implies
$$\nu(\pi_*(T),0)=\lim_{r\rightarrow 0} \frac{1}{2^p\pi^pr^{2p}}
\langle T, 1_{\pi^{-1}(\Ball(0,r))}\Psi\rangle \geq \nu(T,0),$$
which completes the proof.
\end{proof}
The following lemma completes 
the proof of Theorem 3.
\begin{lemma} The set $X$ is $p$-pseudoconcave.
\end{lemma}
\begin{proof} Assume that $X$ is not $p$-pseudoconcave. Then there exist
an open set $U\subset\subset V$ 
and a holomorphic map $f:V'\longrightarrow \C^p$ such that
$f(X\cap U)\not \subset \C^p\setminus \Omega$, where $V'$ is a neighbourhood of 
$\overline U$ and $\Omega$ is 
the unbounded component of $\C^p\setminus f(X\cap bU)$.
\par
Consider $\Phi: V'\longrightarrow \C^p\times V'$ the holomorphic map
given by $\Phi(z):=(f(z),z)$. 
We next choose a domain $U'\subset\subset \C^p\times V'$ such that 
$\Phi(U)$ is a submanifold of $U'$. Then $T':=1_{U'}\Phi_*(T)$ is a 
positive plurisubharmonic
current in $U'$. Moreover it is easy to see that $\nu(T',a)\geq\delta'$ 
for some $\delta'>0$
and for every $a\in X':=\supp(T')$. Let $\pi:\C^{p+n}\longrightarrow \C^p$ be 
the linear projection on the first $p$ coordinates. We have 
$\pi(\overline X'\cap bU')=f(X\cap bU)$ and 
$\pi(X'\cap U')= f(X\cap U)\not \subset \C^p\setminus\Omega$. 
The open set $\Omega$ is
also the unbounded component of $\C^p\setminus 
\pi(\overline X'\cap bU')$.
\par
Therefore $\pi_*(T')$ defines a positive plurisubharmonic
current in $\Omega$ which is null in $\C^p\setminus\pi(\overline X')$.
Hence there is a positive plurisubharmonic
function $\psi$ on $\Omega$, 
null on $\C^p\setminus\pi(\overline X')$
such that $\pi_*(T')=\psi[\Omega]$ in $\Omega$.
Thus for every $x\in \Omega$ we have $\nu(\pi_*(T'),x)=\psi(x)$. 
By Lemma 2, we have 
$\psi(x)\geq \delta'$ for every $x\in \pi(X')\cap\Omega$.
\par
Now fix a point $x$ in the boundary $\Sigma$ of 
$\pi(X')$ in $\Omega$. Let $c\in \Omega$ be a point 
close to $x$. Let $b\in \pi(X')$ such that 
$\dist(b,c)=\dist(\pi(X'),c)$. Since $c$ is 
close to $x$, we have $b\in\Sigma$ and $\Ball(c,|b-c|)\subset
\Omega$. 
By submean
property of plurisubharmonic functions and the fact that $\Psi=0$ on 
$\Ball(c,|b-c|)$,
we obtain 
$$\psi(b) = \lim_{r\rightarrow 0}
\frac{1}{2^p\pi^p r^{2p}}\int_{\Ball(b,r)}\psi(z) \Psi = 
\lim_{r\rightarrow 0}
\frac{1}{2^p\pi^p r^{2p}}\int_{\Ball(b,r)\setminus 
\Ball(c,|b-c|)}\psi(z) \Psi.$$
By uppersemicontinuity, the last limit is smaller or equal to
$$ \lim_{r\rightarrow 0}
\frac{1}{2^p\pi^p r^{2p}}\int_{\Ball(b,r)\setminus \Ball(c,|b-c|)}\psi(b) \Psi
=\frac{\psi(b)}{2}.$$
Thus $\psi(b)\leq \frac{1}{2}\psi(b)$ and $\psi(b)\leq 0$. This contradicts 
the fact that $\varphi(b)\geq\delta'$.
\end{proof}
\begin{corollary} Let $T$ be a positive plurisubharmonic current 
of bidimension $(p,p)$ in a complex manifold $V$ of dimension $n$. 
If $T$ is locally rectifiable, then there exists a locally finite family of  
complex subvarieties $(X_i)_{i\in I}$ 
of pure dimension $p$ of $V$ 
and positive integers $n_i$ such that $T=\sum_{i\in I}n_i [X_i]$. 
In particular, $T$ is closed.  
\end{corollary}
\begin{remark}\rm 
King has proved the same result for rectifiable closed 
positive current \cite{King}. Harvey, Shiffman and
Alexander \cite{Alexander2, HarveyShiffman} proved it for rectifiable 
closed currents (in this case, 
the positivity is not necessary). 
\end{remark}
\begin{proof} By \cite[4.1.28(5)]{Federer}, 
$\nu(T,a)$ is a strictly positive integer for 
$a$ in a dense subset of $X:=\supp(T)$. This, combined with Theorem 3, 
implies that $X$ is a complex 
subvariety of pure dimension $p$ of $V$ and $T=\varphi[X]$, where $\varphi$ is a
weakly plurisubharmonic function on $X$. Since $T$ is 
rectifiable, the function 
$\varphi$ has integer values $\H^{2p}$-almost everywhere. 
Therefore $\varphi$ is essentially 
equal to a positive integer in each 
irreducible component of $X$. 
\end{proof}
The following theorem is a variant of Theorem 3 (\see also Example 1).
\\
\ 
\\
{\bf Theorem 3'}
{\it 
Let $T$ be a positive plurisubharmonic current 
of bidimension $(p,p)$ in a complex manifold $V$ of dimension $n$. 
Assume that $\supp(T)$ has locally finite $\H^{2p}$ measure and the set
$\E:=\{z\in \supp(T),\ \nu(T,z)=0\}$ has zero $\H^{2p-1}$ measure.
Then $\supp(T)$ is a complex subvariety of pure 
dimension $p$ of $V$ and there exists a weakly plurisubharmonic 
function $\varphi$ on $X:=\supp(T)$ such that $T=\varphi[X]$.  
}
\begin{proof} We only need to prove Lemma 3. More precisely, we have to find a point
$b\in \Sigma\setminus \pi(\E)$ such that the upper-density 
$$\overline\Theta(b):=\limsup_{r\rightarrow 0} 
\frac{1}{2^p\pi^p r^{2p}}\int_{\Ball(b,r)\cap\Omega} \Psi$$
of $\Omega$ at $b$ is strictly positive. 
Following \cite[5.8.5, 5.9.5]{Ziemer}, the set 
$\Sigma':=\{b\in\Sigma, \overline\Theta(b)>0\}$ has positive $\H^{2p-1}$ measure. It is
sufficient to take $b$ in $\Sigma'\setminus \pi(\E)$. The last set is not empty since
$\H^{2p-1}(\E)=0$.
\end{proof}
For plurisuperharmonic currents, we have the following theorem.
\begin{theorem} Let $T$ be a positive plurisuperharmonic current 
of bidimension $(p,p)$ in a complex manifold $V$ of dimension $n$. 
Let $K$ be a compact subset of $V$ which admits a Stein neighbourhood.
Assume that in $V\setminus K$ the support
$\supp(T)$ of $T$ has locally finite 
$2p$-dimensional Hausdorff measure.
Then $\supp(T)$ is a complex subvariety of pure 
dimension $p$ of $V$ and there exists a weakly plurisuperharmonic 
function $\varphi$ on $X:=\supp(T)$ such that $T=\varphi[X]$. Moreover, if $T$ is 
pluriharmonic, then $\varphi$ is weakly pluriharmonic on $X$.   
\end{theorem}
\begin{proof} 
By Proposition 1, $X:=\supp(T)$ is $p$-pseudoconcave subset of $V$. 
Moreover, by Corollary 1, $X$ is
a complex subvariety of pure dimension $p$ of $V$. The theorem follows. 
\end{proof}
We remark here that Theorem 4 is false for positive 
plurisubharmonic currents. Consider 
an example.
\begin{example} \rm
Let $\psi$ be a positive subharmonic function in 
$\C$ which is null in the unit disk $\Ball(0,1)$.
We can take for example $\psi(z):=\log|z|$ if $|z|>1$ 
and $\psi(z):=0$ if $|z|\leq 1$.
Let $f$ be a holomorphic function in $\C\setminus \Ball(0,1/2)$ which cannot be
extended to a meromorphic function on $\C$. 
Denote by $Y\subset \C^2$ the graph
of $f$ over $\C\setminus \Ball(0,1/2)$. Let $\varphi$ be the subharmonic 
function on $Y$ given by $\varphi(z):=\psi(z_1)$. Then $T:=\varphi[Y]$ is a 
positive plurisubharmonic current of $\C^2$. Its support $X\subset Y$ is not a
complex subvariety of $\C^2$ and cannot be extended to a complex subvariety.
We can construct more complicate examples by taking countables combinations of 
such currents. 
\end{example}
\begin{corollary} Under the hypothesis of Theorem 3, Theorem 3' or Theorem 4, 
if $V$ is compact,
then there are complex subvarieties $X_1$, $\ldots$, $X_k$ of pure dimension 
$p$ of $V$
and positive real numbers $c_1$, $\ldots$, $c_k$ such that $T=c_1 [X_1] + 
\cdots + c_k[X_k]$. In particular, $T$ is closed.
\end{corollary}
\begin{proof}
By Theorems 3, 3' and 4, $X$ is a complex subvariety 
of pure dimension $p$. Since 
$V$ is compact, $X$ has a finite number of irreducible 
components and $\varphi$ 
is essentially constant in each component.
\end{proof}
%
%\begin{lemma} Let $X$ be a complex subvariety of a 
%compact complex manifold $V$.
%Then every weakly plurisubharmonic function on $X$ 
%is essentially constant in each 
%irreducible component of $X$.
%\end{lemma} 
%
%
\section{Further remarks}
Denote by $\Delta$ the closed unit disk of the complex line.
We say that $V$ satisfies the {\it $p$-Kontinuit\"atssatz} if 
for any smooth family of $p$-dimensional 
closed holomorphic polydiscs $\Delta^p_t$ 
in $V$ indexed by 
$t\in [0,1)$ such that $\cup b\Delta^p_t$ lies on a compact subset of $V$, then 
$\cup \Delta^p_t$ lies on a compact subset of $V$. The following definition 
of $p$-pseudoconcavity does not change the main results of this article.
\\
\ 
\\
{\bf Definition 1'} 
A closed  subset $X$ of a complex manifold $V$ is called {\it $p$-pseudoconcave} if 
$V\setminus X$ satisfies the local $p$-Kontinuit\"atssatz, \ie for every $a\in X$ there is a
neighbourhood $U$ of $a$ such that $U\setminus X$ satisfies the $p$-Kontinuit\"atssatz.
\\
\ 
\\
Let $\Gamma$ be an oriented 
compact real submanifold of dimension $2p-1$ of $\C^n$. Alexander and
Wermer \cite{AlexanderWermer} 
introduced an integer number $\link(\Gamma,H)$ 
(called linking number) 
which depends
on $\Gamma$ and on a complex subvariety $H$ of dimension $n-p$ such that 
$H\cap\Gamma=\emptyset$. One can define the hull of $\Gamma$ by
$$\hull^*(\Gamma,p):=\Gamma\cup\{ z\in \C^n\setminus\Gamma,\ \link(\Gamma,H)\not=0
\mbox{ for every } H \mbox{ containing } a \}.$$
Theorem 1 is valid for this hull but the definition cannot be extended for every 
compact set $\Gamma$. We now try to give an analogous definition.
\par
Let $\Gamma\subset\C$ be a compact set and 
$\Pi:\C^n\longrightarrow \C^p$ be a polynomial map. 
Denote by $\Omega$ the unbounded 
component of $\C^p\setminus\Pi(\Gamma)$ and set 
$\Gamma_\Pi:=\Pi^{-1}(\C^p\setminus\Omega)$.
Then we can define
$$\hull^{**}(\Gamma,p):=\bigcap_\Pi\Gamma_\Pi.$$
Observe that $\link(\Gamma,\Pi^{-1}(\xi))=0$ if $\Gamma$ 
is an oriented compact real 
manifold of dimension $2p-1$ and $\xi\in\Omega$. 
Moreover we have $\hull(\Gamma,p)\subset\hull^{**}(\Gamma,p)$ and 
$\widehat\Gamma = \hull(\Gamma,1)=\hull^{**}(\Gamma,1)$. 
If $\Gamma\cap L$ is a closed irredutible curve with finite length for almost every 
complex $(n-p+1)$-plane $L$, then by slicing method 
we can prove that $\hull(\Gamma,p)=\hull^{**}(\Gamma,p)$. In general, we do not know
when $\hull(\Gamma,p)=\hull^{**}(\Gamma,p)$. There are different problems which 
are related to this question, for example: find topology properties 
of $\hull(\Gamma,p)$,
study $\hull(\Gamma,p)$ when $\Gamma$ is a
compact CR-manifolds, describe $\hull(\Gamma,p)$ using
positive currents...  We also hope that some notions introduced in this article 
can be extended and invertigated in the case of almost complex manifolds.
\par
About positive plurisubharmonic currents, 
the description of the set $\{\nu(T,a)>0\}$ is 
far to be complete. Assume that $\nu(T,a)\geq\delta>0$ for every $a$ 
such that $\nu(T,a)\not=0$; we do not know if $A$ is an analytic set.
\\
\ 
\\
{\bf Acknowledgments.} I would like to thank the Alexander 
von Humboldt foundation 
for its support, Professor J\"urgen Leiterer and Dr. George Marinescu
for their great hospitality and Dr. Nguyen Viet Anh for many help 
during the preparation of this paper. 
Tien-Cuong Dinh\\
Math\'ematique - B\^atiment 425 \\
UMR 8628, Universit\'e Paris-Sud \\
91405 ORSAY Cedex (France) \\
TienCuong.Dinh@math.u-psud.fr
\end{document}